\input amstex
\documentstyle{amsppt}
\magnification=\magstep1 
\overfullrule=0pt
\def\nat{{\Bbb N}}
\def\real{{\Bbb R}}

\def\varep{\varepsilon}
\def\iitem{\itemitem}
\def\trivert{{|\!|\!|}}
\def\ds{\displaystyle}
\def\nto{\buildrel n\over\to}
\def\notnto{\ {\nto\kern-1em\hbox{$\slash\,$}}\ }
\topmatter
\title On impossible extensions of Krivine's Theorem\endtitle
\author E. Odell and Th. Schlumprecht\endauthor
\address Department of Mathematics,
The University of Texas at Austin, 
Austin, TX 78712-1082
\endaddress
\email odell\@math.utexas.edu\endemail 
\address 
Department of Mathematics, 
Texas A\&M University, 
College Station, TX 77843
\endaddress
\email schlump\@math.tamu.edu\endemail 
\thanks Research of E.~Odell supported by NSF and TARP 235. 
Research of Th.~Schlumprecht supported by NSF.\endthanks 
\abstract 
We give examples of two Banach spaces. One Banach space has no spreading 
model which contains $\ell_p$ ($1\le p<\infty$) or $c_0$. The other space 
has an unconditional basis for which $\ell_p$ ($1\le p<\infty$) and $c_0$ 
are block finitely represented in all block bases.
\endabstract
\endtopmatter 

\document
\baselineskip=18pt			

A famous theorem by J.L.~Krivine \cite{K} can be stated as 

\proclaim{Theorem 0.1}
Let $C\ge1$, $n\in\nat$ and $\varep>0$. There exists $m=m (C,n,\varep)\in\nat$ 
so that if $(x_i)_{i=1}^m$ is a finite basic sequence in some Banach space 
with basis constant $C$ then there exist $1\le p\le\infty$ and a block 
basis $(y_i)_1^n$ of $(x_i)_1^m$ so that $(y_i)_1^n$ is 
$(1+\varep)$-equivalent to the unit vector basis of $\ell_p^n$. 
\endproclaim 

Actually this is a stronger version of Krivine's original theorem due to 
Lemberg \cite{L} and H.~Rosenthal \cite{R} (see also \cite{MS} for a 
nice exposition of the proof). Rosenthal also proved 

\proclaim{Theorem 0.2}
Let $(x_i)$ be a basic sequence in a Banach space. There exist a block 
basis $(y_i)$ of $(x_i)$ and a closed set $I\subseteq [1,\infty]$ such 
that if $p\in I$ and $(z_i)$ is any block basis of $(y_i)$, then $\ell_p$ 
is block finitely represented in $(z_i)$.
\endproclaim 

He raised 

\proclaim{Problem 0.3} 
Does $I= \{p\}$ for some $p$?
\endproclaim 

We show in \S2 that this is not the case. 
In our example $I= [1,\infty]$. In fact we construct an unconditional 
basic sequence $(x_i)$ with the property that every 1-unconditional basic 
sequence is block finitely represented in every block basis of $(x_i)$. 

The second problem we address involves spreading models. Not every 
infinite dimensional Banach space must contain $c_0$ or $\ell_p$ for some 
$1\le p<\infty$ as shown by Tsirelson \cite{T}. 
Krivine's theorem gives certain finite information about basic sequences. 
Between these two results lies the well known 

\proclaim{Problem 0.4} 
Let $X$ be an infinite dimensional Banach space. Does $X$ have $c_0$ 
or $\ell_p$ (for some $1\le p<\infty$) as a spreading model?
\endproclaim 

In \S1 we exhibit a space $X$ with an unconditional subsymmetric basis 
having the property that if $Y$ is any spreading model of $X$ then $Y$ 
does not contain $c_0$ or $\ell_p$ ($1\le p<\infty$). 

The original space of Tsirelson has $c_0$ as a spreading model and its dual  
space $T$ as described by Figiel and Johnson \cite{FJ} has $\ell_1$ as 
spreading model. Numerous relatives of $T$ have subsequently been defined 
(see \cite{CS}) using variants of the clever implicit description of the 
norm due to Figiel and Johnson but fail to be a counterexample to Problem~0.4. 
The space $S$ \cite{S1,2} comes close but was shown by Pei-Kee Lin to have 
$\ell_1$ as a spreading model (it is not known if $S^*$ has $c_0$ as a 
spreading model). 

Both of our examples are Tsirelson type spaces --- spaces defined by an 
implicit Figiel-Johnson type norm equation --- and involve modifying the 
norm of $S$. The example in \S1 modifies $S$ along the lines of 
W.T.~Gowers' recent example \cite{G}. It is unknown whether Gowers' space 
has $\ell_1$ as a spreading model. 

The theory of spreading models, which originated with the work of Brunel and 
Sucheston \cite{BS1,2}, is now fairly well established. For background 
information see \cite{BL} (or \cite{O} for a quick introduction). 

Our terminology is standard as may be found in \cite{LT}. 
If $A\subseteq X$ where $X$ is 
a Banach space then span~$(A)$ is the linear span of $A$. 
$S(X)$ is the unit sphere of $X$ and $Ba(X)$ is the unit ball of $X$. 
A basic sequence $(x_i)$ is {\it block finitely represented\/} in $(y_i)$ 
if for all $\varep>0$ and $n\in\nat$ there exists a block basis 
$(z_i)_{i=1}^n$ of $(y_i)$ satisfying 
$$(1+\varep)^{-1}\Big\| \sum_1^n a_ix_i\Big\| 
\le \Big\| \sum_1^n a_iz_i\Big\| 
\le (1+\varep) \Big\|\sum_1^n a_ix_i\Big\|$$ 
for all $(a_i)_1^n\subseteq \real$. 
$\ell_p$ {\it is block finitely represented in\/} 
$(y_i)$ if the unit vector basis 
of $\ell_p$ is block finitely represented in $(y_i)$. 

\head \S1. A space with no spreading model containing $c_0$ or $\ell_p$ 
\endhead 

Let $c_{00}$ be the linear space of all finitely supported real valued 
functions on $\nat$. Let $f(i) = \log_2 (1+i)$ for $i\in\nat$. 
For $E,F\subseteq\nat$ we write $E<F$ if $\max E< \min F$. For $x\in c_{00}$ 
and $E\subseteq \nat$, let $Ex\in c_{00}$ be given by $Ex(i) = x(i)$ 
if $i\in E$ and $0$ otherwise. Fix an increasing sequence of integers 
$(n_k)$ with 
$$\sum_{k=1}^\infty {1\over f(n_k)}  < \tfrac1{10}\ .
\leqno(1.0)$$ 

\proclaim{Proposition 1.1} 
There exists a $1$-unconditional norm $\|\cdot\|$ on $c_{00}$ which 
satisfies the implicit equation 
$$\leqalignno{
&\|x\| = \max \biggl\{ \|x\|_{c_0},\biggl( \sum_{k=1}^\infty \|x\|_{n_k}^2 
\biggr)^{1/2}\biggr\}\ \text{ where}
&(1.1)\cr 
&\|x\|_k = \max_{E_1<\cdots < E_k} {1\over f(k)} \sum_{i=1}^k \|E_ix\|\ .
&(1.2)\cr}$$ 
\endproclaim 

\demo{Proof}
We follow the standard Tsirelson norm construction of \cite{FJ}. 
Let $\|x\|_{(0)} \equiv \|x\|_{c_0}$. 
If $\|x\|_{(k)}$ has been defined set 
$$\|x\|_{(k+1)} = \max \biggl\{ \|x\|_{(k)}, \biggl( \sum_{i=1}^\infty 
\|x\|_{(k,n_i)}^2\biggr)^{1/2}\biggr\}$$ 
where 
$$\|x\|_{(k,i)} = \max_{E_1<\cdots < E_i} {1\over f(i)} 
\sum_{j=1}^i \|E_jx\|_{(k)}$$ 
$\|\cdot\|_{(k)}$ is a norm for each $k$ with $\|x\|_{(0)} \le 
\|x\|_{(1)}\le\cdots$ and all norms are dominated by 
$\|\cdot\|_{\ell_1}$. 
The latter fact can be seen from 
observing that $\|e_j\|_{(k)} = \|e_j\|_{(0)} =1$ for all $k$, where  
$(e_i)$ is the unit vector basis for $c_{00}$. 
The proposition follows by taking $\|x\|\equiv \lim_k \|x\|_{(k)}$.\qed 
\enddemo 

We let $X$ be the Banach space given by completing the space of 
Proposition~1.1. The unit vector basis $(e_i)$ is a normalized 1-unconditional 
subsymmetric basis for $X$. 

\proclaim{Proposition 1.2} 
Let $(x_i)$ be a normalized block basis of $(e_i)$ with spreading model 
$(u_i)$. Let $U= [(u_i)]$. Then $\ell_p$ for 
$1<p<\infty$ and $c_0$ are not block 
finitely representable in $(u_i)$.
\endproclaim 

\demo{Proof} 
Let $1<p<\infty$ (a similar argument works for $c_0$). 
If $\ell_p$ is block finitely representable in $(u_i)$, then 
$\ell_p$ is block finitely representable in $(x_i)$. 
But if $(y_j)_{j=1}^{n_i}$ is a normalized block basis of $(e_i)$, 
then by (1.1) 
$$\Big\| \sum_{j=1}^{n_i} y_j\Big\| 
\ge \Big\| \sum_{j=1}^{n_i} y_j\Big\|_{n_i}  
\ge {1\over f(n_i)} n_i$$ 
which shows this to be impossible.\qed 
\enddemo 

Thus by Krivine's theorem we need only show that such a $U$ cannot 
contain $\ell_1$. 

\proclaim{Theorem 1.3} 
Let $(x_i)$ be a normalized block basis of $(e_i)$ with spreading model 
$(u_i)$. Then $U\equiv [(u_i)]$ does not contain $\ell_1$. 
\vskip1pt
In particular $X$ cannot contain $\ell_1$ or $c_0$ and thus is reflexive. 
\endproclaim 

\demo{Proof} 
We first prove that $(u_i)$ cannot be equivalent to the unit vector basis of 
$\ell_1$. If not then we may assume (by replacing $(x_i)$ by a suitable 
bounded length block basis) that $\|x_i\|_{c_0}<1$ for all $i$ and for 
all $(c_i)\subseteq\real$, 
$$\big\|\sum c_i u_i\big\| \ge (.99) \sum |c_i|\ .$$ 
This follows from James' proof that $\ell_1$ is not distortable 
(see e.g., \cite{BL}, p.43). 

For $i,k\in\nat$ let $d_{i,k} = \|x_i\|_{n_k}$. 
Then $d_i\equiv (d_{i,k})_{k=1}^\infty \in S(\ell_2)$, the unit sphere 
of $\ell_2$, for all $i\in\nat$. By passing to a subsequence 
of $(x_i)$ we may assume 
that $(d_i)$ converges weakly to $d\equiv (a_i) \in Ba(\ell_2)$ and 
$$\lim_{\scriptstyle \ell\to\infty\atop\scriptstyle \ell<i_1<\cdots < i_n} 
\Big\| \sum_{j=1}^n d_{i_j}\Big\|_{\ell_2}\ \text{ exists for all }\ n\ .$$
Fix an integer $n\ge2$. Then 
$$\eqalign{.99n \le \Big\| \sum_{i=1}^n u_i\Big\| 
&=\lim_{\scriptstyle \ell\to\infty\atop\scriptstyle \ell<i_1<\cdots< i_n} 
\Big\|\sum_{j=1}^n x_{i_j}\Big\|\cr
&=\lim_{\scriptstyle \ell\to\infty\atop \scriptstyle\ell<i_1<\cdots<i_n} 
\biggl( \sum_{k=1}^\infty \Big\|\sum_{j=1}^n x_{i_j}\Big\|_{n_k}^2
\biggr)^{1/2}\cr 
&\le\lim_{\scriptstyle\ell\to\infty\atop\scriptstyle \ell<i_1<\cdots<i_n} 
\Biggl( \sum_{k=1}^\infty \biggl(\sum_{j=1}^n d_{i_j,k}\biggr)^2 
\Biggr)^{1/2}\cr 
&= \lim_{\scriptstyle\ell\to\infty\atop\scriptstyle \ell<i_1<\cdots<i_n} 
\Big\| \sum_{j=1}^n d_{i_j} \Big\|_{\ell_2}\ .\cr}$$ 
Let $\|d\|^2 +\varep^2 =1$. Since $(d_i) \subseteq S(\ell_2)$ converges 
weakly to $d$, the latter limit above is 
$$= (n^2 \|d\|^2 + n\varep^2)^{1/2}\ .$$ 
This can be seen by the standard gliding hump argument, choosing $d_{i_j}$ 
essentially equal to $d +h_j$ where $\|h_j\|=\varep$ and the ``humps'' 
$(h_j)$ are disjointly supported in $\ell_2$ and nearly disjoint from $d$. 
Thus $.99 \le (\|d\|^2+\varep^2/n)^{1/2}$. Since $n$ was arbitrary, 

\iitem{i)} $\|d\|\ge .99$ and hence since $\|d\|^2 + \varep^2=1$, 
$\varep < .15$. 

\noindent Fix $k_0\in \nat$ so that 
\smallskip
\iitem{ii)} $\ds\biggl(\sum_{k=1}^{k_0} a_k^2\biggr)^{1/2} > .98$ and hence 
\smallskip
\iitem{iii)} $\ds\biggl(\sum_{k=k_0+1}^\infty a_k^2\biggr)^{1/2} \le .2$ 
\smallskip
\noindent Since $(d_{i,k})_{k=k_0+1}^\infty$ converges weakly to 
$(a_k)_{k=k_0+1}^\infty$ in $\ell_2$ and $\lim_{i\to\infty} 
\|(d_{i,k})_{k= k_0+1}^\infty - (a_k)_{k=k_0+1}^\infty\| \break = \varep$, 
the analysis above yields (using iii)) that, 
$$\displaylines{
\lim_{\scriptstyle\ell\to\infty\atop\scriptstyle \ell<i_1<\cdots<i_n} 
\biggl( \sum_{k=k_0+1}^\infty \Big\| \sum_{j=1}^n x_{i_j} \Big\|_{n_k}^2 
\biggr)^{1/2} 
\le \big( n^2 (.2)^2 + n\varep^2 \big)^{1/2}\cr 
< .2\big(n^2 +n\big)^{1/2} < .3n\ \text{ if $n$ is sufficiently large.}
\cr }$$
Choose and fix $N$ so that this holds and 

\iitem{iv)} $2n_{k_0}^2 < .1 N^2$. 

\noindent Let $1>\delta >0$ be such that 
\smallskip
\iitem{v)} $\ds \sum_{k=1}^{k_0} (a_k+\delta)^2 <2\ .$
\smallskip
\noindent Then choose $i_1<i_2<\cdots<i_N$ so that 
\smallskip
\iitem{vi)} $\ds\biggl(\sum_{k=k_0+1}^\infty \Big\|\sum_{j=1}^N x_{i_j} 
\Big\|_{n_k}^2\biggr)^{1/2} < .3N\ ,$ 
\smallskip 
\iitem{vii)} $\big|\ \|x_{i_j}\|_{n_k} -a_k\big| <\delta$ for $j\le N$, 
$k\le k_0$, and 
\smallskip
\iitem{viii)} $\ds \Big\| \sum_{j=1}^N x_{i_j}\Big\| 
= \biggl( \sum_{k=1}^\infty \Big\|\sum_{j=1}^N x_{i_j}\Big\|_{n_k}^2\biggr)^
{1/2} \ge .98N$ 
\smallskip
\noindent From viii) and vi) we obtain that 
\smallskip
\iitem{ix)} $\ds\biggl(\sum_{k=1}^{k_0}
\Big\|\sum_{j=1}^N x_{i_j}\Big\|_{n_k}^2 \biggr)^{1/2}\ge .9N$. 

Let $k\le k_0$ be fixed. Then 
$$\Big\|\sum_{j=1}^N x_{i_j}\Big\|_{n_k} 
= {1\over f(n_k)} \sum_{\ell=1}^{n_k} \Big\| E_\ell \biggl( \sum_{j=1}^N 
x_{i_j}\biggr)\Big\|$$ 
for some choice $E_1< E_2 <\cdots < E_{n_k}$. 
Now if 
$$I= \bigl\{ j\le N :E_\ell (x_{i_j})\ne0\text{ for at least two different  
$\ell$'s }\bigr\}\ ,$$ 
then $|I| \le n_k$. Thus 
$$\eqalign{ 
{1\over f(n_k)} \sum_{\ell=1}^{n_k} \Big\| E_\ell\biggl( \sum_{j=1}^N 
x_{i_j}\biggr)\Big\| 
&\le {1\over f(n_k)} \Biggl( \sum_{\ell=1}^{n_k} \biggl( \sum_{j\in I} 
\|E_\ell x_{i_j}\| + \Big\| E_\ell\biggl( \sum_{\scriptstyle j=1\atop 
\scriptstyle j\notin I}^N x_{i_j}\biggr)\Big\|\biggr)\Biggr)\cr 
&\le \sum_{j\in I} {1\over f(n_k)} \sum_{\ell=1}^{n_k} 
\|E_\ell x_{i_j}\| + {1\over f(n_k)} \sum_{\scriptstyle j=1\atop 
\scriptstyle j\notin I}^N \|x_{i_j}\|\cr 
&\le \sum_{j\in I} \|x_{i_j}\|_{n_k} + {N\over f(n_k)}\cr 
&\le |I| (a_k +\delta) + {N\over f(n_k)}\qquad\text{(by vii))}\cr 
&\le n_k (a_k+\delta) +{N\over f(n_k)}\ .\cr}$$ 
We thus obtain, 
$$\eqalign{ 
\sum_{k=1}^{k_0} \Big\|\sum_{j=1}^N x_{i_j}\Big\|_{n_k}^2 
&\le \sum_{k=1}^{k_0} \left( n_k (a_k+\delta) +{N\over f(n_k)}\right)^2\cr 
&\le  n_{k_0}^2 \sum_{k=1}^{k_0} (a_k+\delta)^2 + 2n_{k_0} N 
\sum_{k=1}^{k_0} {(a_k+\delta) \over f(n_k)} 
+N^2\sum_{k=1}^{k_0} {1\over f(n_k)^2}\cr 
&< 2n_{k_0}^2 + 4n_{k_0} N(.1) + .1 N^2\cr}$$ 
by v), $|a_k+\delta|\le 2$ and (1.0). 
This in turn is by iv), 
$$< 2n_{k_0}^2 + .5N^2 < .6N^2\ .$$ 
Thus 
$$\biggl( \sum_{k=1}^{k_0} \Big\| \sum_{j=1}^N x_{i_j}\Big\|_{n_k}^2
\biggr)^{1/2} < .8N$$ 
which contradicts ix). 

We next show that $U$ does not contain $\ell_1$. 
By a diagonal argument and passing to a subsequence of $(x_i)$ if necessary 
we may assume that $(x_i)$ has a spreading model for each norm 
$\|\cdot\|_{n_k}$; i.e., 
$$\lim_{\scriptstyle\ell\to\infty\atop 
\scriptstyle \ell<i_1<\cdots<i_m} \|\sum_{j=1}^m b_j x_{i_j}\|_{n_k} 
\text{ exists for all }(b_j)_1^m \subseteq \real\ .$$ 
Let $v_i = \sum_{j=m_i+1}^{m_{i+1}} b_j u_j$ be a normalized block basis 
of $(u_i)$ with $\|\sum c_i v_i\|\ge .99\sum |c_i|$ for $(c_i)\subseteq 
\real$. Let $d_i = (d_{i,k})_{k=1}^\infty \in S(\ell_2)$ be defined by 
$$d_{i,k} = \lim_{\scriptstyle \ell\to\infty\atop\scriptstyle \ell < 
i_{m_i+1} < \cdots < i_{m_{i+1}}} \Big\| \sum_{j=m_i+1}^{m_{i+1}} 
b_ix_{i_j}\Big\|_{n_k}\ .$$ 
We may assume by passing to a subsequence of $(v_i)$ that 
$(d_i)$ converges weakly to $d\in Ba (\ell_2)$. 
In fact the entire argument above  now carries over; $\|d\|\ge .99$ and 
$d$ determines $k_0$ as before. This then defines $N$ and $\delta$ 
and yields a contradiction.\qed
\enddemo

Our work above suggests two natural problems. 
Let us say $E\to F$ if $F$ is a spreading model of some basic sequence 
in $E$ and $E\nto F$ if $E\to E_1\to \cdots\to E_n\to F$ 
for some sequence of Banach spaces $E_1,\ldots,E_n$. 

\proclaim{Problem 1.4} 
Given $n\ge2$ does there exist a Banach space $E$ such that if $F=\ell_p$ 
$(1\le p<\infty)$ or $c_0$ then $E\notnto F$? 
\endproclaim 

\proclaim{Problem 1.5} 
Does there exist a Banach space $E$ such that for all $n$, 
$E\notnto F$ whenever $F=\ell_p$ $(1\le p<\infty)$ or $c_0$?
\endproclaim 

\head \S2 \endhead 

In this section we prove 

\proclaim{Theorem 2.1} 
There exists a $1$-unconditional basic sequence $(e_i)$ such that if 
$n\in\nat$ and $(x_i)_{i=1}^n$ is a finite $1$-unconditional basic 
sequence, $\varep>0$ and $(y_i)_{i=1}^\infty$ is a block basis of 
$(e_i)$ then there exists a finite block basis $(z_i)_{i=1}^n$ of $(y_i)$ 
which is $(1+\varep)$-equivalent to $(x_i)_{i=1}^n$.
\endproclaim 

We first observe that it is not necessary to directly check all such 
sequences $(x_i)_{i=1}^n$. Let $n\in\nat$ and let $B(\ell_\infty^n)$ 
denote the $n^2$-dimensional Banach space of all bounded linear operators 
on $\ell_\infty^n$. $B(\ell_\infty^n)$ has a matrix basis 
$(e_{i,j})_{i,j=1}^n$ given by $e_{i,j}(f_k)=f_j$ for $i=k$ and $0$ otherwise 
(where $(f_k)_{k=1}^n$ is the unit vector basis of $\ell_\infty^n$). 

It is routine to check that 
$$\Big\| \sum_{i,j=1}^n a_{i,j} e_{i,j}\Big\|  
= \max_{j\le n} \sum_{i=1}^n |a_{i,j}| 
\leqno(2.1)$$ 
and thus $(e_{i,j})_{i,j=1}^n$ is a $1$-unconditional basis for 
$B(\ell_\infty^n)$. We order this basis lexicographically: 
$(e_{11}, e_{12},\ldots,e_{1n},e_{21}, e_{22},\ldots,e_{nn})$. 

\proclaim{Proposition 2.2} 
Let $(y_k)_{k=1}^m$ be a $1$-unconditional basic sequence for some $m\in\nat$ 
and let $\varep>0$. There exists $n\in\nat$ and a block basis $(x_k)_{k=1}^m$ 
of the basis $(e_{i,j})$ for $B(\ell_\infty^n)$ so that $(x_k)_{k=1}^m$ 
is $(1+\varep)$-equivalent to $(y_k)_{k=1}^m$.
\endproclaim 

\demo{Proof} 
Without loss of generality we may assume that span$_{1\le i\le m} (y_i)$ 
is a subspace of $\ell_\infty^n$ for some $n\in\nat$. 
Write $y_k = \sum_{j=1}^n a_{k,j} f_j$ for $k\le m$ and define 
$$x_k = \sum_{j=1}^n a_{k,j} e_{k,j}\  .$$ 
Let $(b_k)_{k=1}^m$ be scalars. Then from (2.1), 
$$\eqalign{\Big\| \sum_{k=1}^m b_k x_k\Big\| 
& = \Big\| \sum_{k=1}^m \ \sum_{j=1}^n b_k a_{k,j} e_{k,j}\Big\|\cr 
& = \max_{j\le n} \sum_{k=1}^m |b_k a_{k,j}|\ .\cr}$$ 
Also using the 1-unconditionality of $(y_k)$, 
$$\eqalignno{\Big\| \sum_{k=1}^m b_k y_k\Big\| 
&= \Big\|\sum_{k=1}^m\ \sum_{j=1}^n b_k a_{k,j} f_j\Big\|\cr 
&= \max_{j\le n} \sum_{k=1}^m |b_k a_{k,j}|\ .&\qed\cr}$$
\enddemo 

Before proceeding we set some notation. 
For $t>0$, set, as in \S1, $f(t) = \log_2 (t+1)$. 
A finite sequence of pairs 
$((m_i,E_i))_{i=1}^k$, where $m_1<m_2<\cdots < m_k$ are integers in $\nat$ 
and $E_1<\cdots < E_k$ are finite subsets of $\nat$, is {\it admissible\/} 
if $m_1\ge2$ and $f (m_{i+1}) >\sum_{j=1}^i |E_j|$ for $1\le i<k$. 

\proclaim{Proposition 2.3} 
There is a norm $\|\cdot\|$ on $c_{00}$ satisfying the following 
implicit equation. 
$$\|x\| = \max \Biggl\{ \|x\|_\infty ,\sup \biggl\{ {1\over f(\ell)} 
\sum_{i=1}^\ell \trivert E_ix\trivert_{m_i} : \ell\in \nat\ ,\ 
(m_i,E_i)_{i=1}^\ell\text{ is admissible}\biggr\}\Biggr\}
\leqno(2.2)$$ 
where for $m\ge2$, $\trivert\cdot\trivert_m$  is a norm on $c_{00}$ 
given by 
$$\trivert x\trivert_m = \sup \biggl\{ {1\over m} \sum_{i=1}^m 
\|F_ix\| : F_1 < \cdots < F_m\biggr\}\ .
\leqno(2.3)$$ 
\endproclaim 

\demo{Proof} 
The proof is similar to the Figiel-Johnson construction of the Tsirelson 
norm. We first inductively define for every $n\in\nat\cup \{0\}$ a norm 
$|\cdot|_n$ on $c_{00}$. Set $|x|_0 = \|x\|_\infty = \max_{i\in\nat} 
|x(i)|$. If $|\cdot|_n$ has been defined, given $m\in\nat$ and $x\in c_{00}$, 
set 
$$|x|_{(n,m)} = \sup \biggl\{ {1\over m} \sum_{i=1}^m |F_ix|_n : 
F_1<\cdots <F_m\}\ .$$ 
Then set 
$$|x|_{n+1} = \max \Biggl\{ |x|_n ,\sup \biggl\{ {1\over f(\ell)} 
\sum_{i=1}^\ell |E_ix|_{(n,m)} : (m_i,E_i)_{i=1}^\ell\text{ is admissible}
\biggr\}\Biggr\}\ .$$ 
Finally set 
$\|x\| = \max_k |x|_k$ and define $\trivert x\trivert_m$ by (2.3). 

We check that this norm satisfies (2.2). Let $x\in c_{00}$ and let 
$((m_i,E_i))_{i=1}^\ell$ be admissible. For $1\le i\le \ell$ let 
$F_1^i <F_2^i <\cdots < F_{m_i}^i$ be subsets of $E_i$. Then 
$$\eqalign{{1\over f(\ell)} \sum_{i=1}^\ell {1\over m_i} 
\sum_{j=1}^{m_i} \|F_j^i x\| 
& = \max_{k\ge0} {1\over f(\ell)} \sum_{i=1}^\ell {1\over m_i} 
\sum_{j=1}^{m_i} |F_j^i x|_k\cr
&\le \max_{k\ge0} {1\over f(\ell)} \sum_{i=1}^\ell |E_ix|_{(k,m_i)}\cr 
&\le \max_{k\ge0} |x|_k = \|x\|\ .\cr}$$ 
Thus $\|x\|\ge$ right side of (2.2). 
If $\|x\|= |x|_0 = \|x\|_\infty$ we have equality. 
Otherwise $\|x\| = |x|_k > |x|_{k-1}$ for some $k\ge1$. 
Thus there exists an admissible collection $((m_i,E_i))_{i=1}^\ell$ so that 
$$\eqalignno{\|x\|  = |x|_k &= {1\over f(\ell)} 
\sum_{i=1}^\ell |E_ix|_{(k-1,m_i)} \cr 
&\le {1\over f(\ell)} \sum_{i=1}^\ell \trivert E_ix\trivert_{m_i} 
\le \text{ right side of (2.2).}
&\qed\cr}$$
\enddemo 

Let $X$ be the completion of $c_{00}$ under the norm of (2.2). 
The unit vector basis $(e_i)$ for $c_{00}$ is a normalized 1-unconditional,  
1-subsymmetric basis for $X$. Thus 
$$\|\sum a_i e_i\| = \|\sum \varep_i a_i e_{k_i}\|$$ 
whenever $\sum a_i e_i\in X$, $\varep_i=\pm 1$ and $k_1<k_2 <\cdots$. 
Furthermore (2.2) and (2.3) hold for all $x\in X$.

The proof of Theorem 2.1 is quite technical and thus we first sketch 
the proof and give the intuition behind the argument. 

It is not difficult to show that each block basis $(y_i)$ of $(e_i)$ 
which has the unit vector basis of $c_0$ block finitely represented must 
also have the unit vector basis of $\ell_1$ block finitely represented 
(see the proof of step~1, below). 
Thus we deduce from this and Krivine's theorem that the unit vector basis 
of $\ell_p$ is block finitely represented in $(y_i)$ for some $1\le p <
\infty$ (step~1). 
We then choose, given $m\in\nat$ and $\varep>0$, $(x_i)_{i=1}^m$ 
where $x_i = n_i^{-1/p}\sum_{j=1}^{n_i} x_{i,j}$ is a rapidly increasing 
sequence of $\ell_p$-averages in $(y_i)$. 
This means that $(x_{11}, x_{12},\ldots,x_{1n_1}, x_{21},\ldots, 
x_{m,n_m})$ is a normalized block basis of $(y_i)$, $(x_{i,j})_{j=1}^{n_i}$ is 
$1+\varep$-equivalent to the unit vector basis of $\ell_p^{n_i}$ and the 
$n_i$'s are rapidly increasing ($n_1$ is large and, $n_i = n_i(\sum_{j=1}^{i-1}
x_j,p,\varep)$).
We deduce from step~4 below that $(x_i)_{i=1}^m$ is 
$(1+3\varep)(1+\varep)$-equivalent 
to the unit vector basis of $\ell_\infty^m$. 

By the observation mentioned at the beginning of this sketch we deduce also 
that $\ell_1$ is block finitely represented in $(y_i)$. Thus we can find,  
given $n$, $\varep$, a block basis $(x(i,j))_{i,j\le n}$ of $(y_i)$ where 
$x(i,j) = {1\over k_0} \sum_{s=1}^{k_0} x(i,j,s)$ and $(x(i,j,s))_{i,j,s}$ 
is a judiciously chosen block basis of $(y_i)$ with 
$(x(i,j,s))_{i=1}^n{}_{,\,s=1}^{,\,k_0}$ $(1+\varep)$-equivalent to the unit 
vector basis of $\ell_1^{nk_0}$ for $j\le n$. 
Finally we show $(x(i,j))_{i,j=1}^n$ is $(1+\varep)$-equivalent to the 
basis $(e_{i,j})_{i,j=1}^n$ given by (2.1). 

We first set some notation. Following \cite{GM}, for $1\le p\le\infty$, 
$n\in\nat$ and $C\ge1$ we call $x$ an {\it $\ell_p^n$-average with constant\/} 
$C$ if $x=n^{-1/p}\sum_{i=1}^n x_i$ where $(x_i)_{i=1}^n$ is a normalized block 
basis of $(e_i)$ which is $C$-equivalent to the unit vector basis of 
$\ell_p^n$. Note that then $C^{-1}\le \|x\|\le C$. 
For $\ell,m_0\in\nat$ and $x\in X$, define 
$$\leqalignno{\qquad 
\|x\|_{(\ell,m_0)} &= \sup \biggl\{ {1\over f(\ell)} \sum_{i=1}^\ell 
\trivert E_ix\trivert_{m_i} :m_1\ge m_0\text{ and } 
(m_i,E_i)_{i=1}^\ell\text{ is admissible}\biggr\}\text{ and}
&(2.5)\cr
\|x\|_\ell &=\sup \biggl\{ {1\over f(\ell)} \sum_{i=1}^\ell 
\trivert E_ix\trivert_{m_i} : (m_i,E_i)_{i=1}^\ell \text{ is admissible}
\biggr\} \equiv \|x\|_{(\ell,2)}
&(2.6)\cr}$$

\demo{Remark 2.4}
Thus $\|x\| = \max \{\|x\|_\infty, \sup_\ell \|x\|_\ell\}$. 
If $x\in c_{00}$ with $\|x\| \ne \|x\|_\infty$, then there exists $\ell>1$ 
so that $\|x\| = \|x\|_\ell$. Indeed suppose $1=\|x\|=\|x\|_1$. 
Then $\|x\| = \trivert x\trivert_m$ for some  $m\ge2$. 
Choose $m$ maximal with this property. 
It follows that $x= \sum_{i=1}^m x_i$ where $(x_i)_1^m$ is a block basis of 
$(e_i)$ and $\|x_i\| = \|x\|=1$ for all $i\le m$. 
Since $m$ was maximal, $\trivert x_1\trivert_k<1$ for all $k$ and so 
$\|x_1\|_\ell=1$ for some $\ell>1$. 
\enddemo 

\proclaim{Lemma 2.5}
Let $n\in\nat$ and let $(x_i)_{i=1}^n$ be a block basis of $(e_i)$ so that 
for each $i\le n$, $x_i$ is an $\ell_1^{k_i}$-average with constant $2$ 
for some $k_1,\ldots,k_n\in\nat$. Let $k_0 = \min \{k_i:1\le i\le n\}$. 
Then for all $\ell\in\nat$ and $(a_i)_{i=1}^n \subseteq [-1,1]$, 
$$\Big\| \sum_1^n a_i x_i\Big\|_\ell 
\le {1\over f(\ell)} \biggl[ \Big\|\sum_{i=1}^n a_i x_i\Big\| 
+ 6\ell nk_0^{-1/2}\biggr]$$ 
\endproclaim 

\demo{Remark 2.6} 
An easy computation shows that Lemma 2.5 implies that if ${f(\ell)-1\over\ell} 
> 12nk_0^{-1/2}$ then if $\|\sum_1^n a_ix_i\|=1$, 
$$\Big\| \sum_1^n a_i x_i\Big\|_\ell < 
\Big\| \sum_1^n a_ix_i\Big\|\ .$$ 
\enddemo 

We postpone the proofs of Lemma 2.5 and the next lemma. 

\proclaim{Lemma 2.7}
Let $L_0\in\nat$ and $1>\varep >0$. 
There exits $L_1,L'_1\in\nat$ with $L_0<L_1<L'_1$ so that for any $m_0\in\nat$ 
and any block basis $(y_i)$ of $(e_i)$ there is an $x\in\text{\rm span}
(y_i)$ satisfying: $\|x\|=1$ and 
\vskip1pt
\iitem{a)} $\|x\|_\ell \le {2\over f(\ell)}$ if $\ell\le L_0$.
\vskip1pt 
\iitem{b)} $\|x\|_{(L_1,m_0)} \ge 1-\varep$. 
\vskip1pt 
\iitem{c)} $\|x\|_\ell \le\varep$ if $\ell\ge L'_1$.
\endproclaim 

\demo{Proof of Theorem 2.1} 
Let $(y_i)$ be a block basis of $(e_i)$, $n\in\nat$ and $1>\varep>0$. 
By (2.1) and Proposition~2.2 it suffices to produce a block basis 
$(x(i,j))_{i,j=1}^n$ of $(y_i)$ (ordered lexicographically) so that 
for all $(a_{i,j})_{i,j=1}^n\subseteq \real$, 
$$(1+\varep)^{-1}\max_{j\le n} \sum_{i=1}^n |a_{i,j}| 
\le \Big\| \sum_{i,j=1}^n a_{i,j} x(i,j)\Big\| 
\le (1+\varep) \max_{j\le n} \sum_{i=1}^n |a_{i,j}|\ . 
\leqno(2.7)$$

Choose $\delta >0$ so that 
$$\leqalignno{
&(1+\varep)^{-1} < (1-\delta)^2\ \text{ and} &(2.8)\cr
&\delta < 3^{-1} \varep n^{-2}\ .&(2.9)\cr}$$
Let $L_0\in\nat$ so that 
$$f(L_0) > \delta^{-1} 
\leqno(2.10)$$ 
and choose $k_0\in\nat$, $k_0> \max (L_0,n)$ so that 
$$k_0^{-1/2} < {f(L_0) -1\over 12nL_0}\ . 
\leqno(2.11)$$ 
Choose $L'_0 \in\nat$, $L'_0 >k_0$ so that 
$${f(L'_0) \over  f(nk_0 L'_0)} > 1-\delta\ . 
\leqno(2.12)$$ 

We then choose $L'_0 < L_1 < L'_1 < L_2 < L'_2 <\cdots < L_n < L'_n$ 
as follows. $L_1 $ and $L'_1$ are chosen as in Lemma~2.7 for 
$L_0$ (of Lemma~2.7) $\equiv L'_0$ and $\varep\equiv\delta$. 
If $L_r < L'_r$ are chosen, choose $L'_r < L_{r+1} < L'_{r+1}$ by 
Lemma~2.7 with $L_0\equiv L'_r$ and $\varep\equiv \delta$. 

Choose now inductively, using Lemma 2.7, a lexicographically ordered 
normalized block basis $\{x(i,j,s): 1\le i,j\le n$, $s\le k_0\}$ of 
$(y_i)$  along with integers $\{m_0 (i,j,s): i,j\le n$, $s\le k_0\}$ 
and an admissible family $\{ m_t^{(i,j,s)} E_t^{(i,j,s)} : 
i,j\le n$, $s\le k_0$, $1\le t\le L_j\}$ so that 
$$\leqalignno{
&m_0(1,1,1) =1\ \text{ and }\ f\bigl(m_0(i,j,s)\bigr) 
> \sum_{(i',j',s') < (i,j,s)} 
\sum_{t=1}^{L'_j} |E_t^{(i',j',s')}|\ .
&(2.13)\cr
&E_t^{(i,j,s)} < E_{t'}^{(i',j',s')} \ \text{ if }\ 
(i,j,s,t) < (i',j',s',t')\ .
&(2.14)\cr 
&\|x (i,j,s)\|_\ell \le {2\over f(\ell)} \ \text{ if }\ \ell\le L'_{j-1}\ .
&(2.15)\cr 
&\|x(i,j,s)\|_{(L_j,m_0(i,j,s))} 
= {1\over f(L_j)} \sum_{t=1}^{L_j} \trivert E_t^{(i,j,s)} x(i,j,s)\trivert_
{m_t^{(i,j,s)}} \ge 1-\delta \ . 
&(2.16)\cr 
&\|x(i,j,s)\|_\ell \le \delta \ \text{ if }\ \ell \ge L'_j\ . 
&(2.17)\cr}$$ 

The choice of $L_1,L'_1$ permits us to choose 
$x(1,1,1)$ satisfying (2.15)--(2.17) for $(i,j,s)=(1,1,1)$. 
Assuming that $x(i,j,s)$ has been chosen for all $(i,j,s)<(i_0,j_0,s_0) 
\le (n,n,k_0)$, let $m_0 (i_0,j_0,s_0)$ be chosen as in (2.13) and choose 
$x(i_0,j_0,s_0)$, and an admissible family $(m_t^{(i_0,j_0,t_0)}, 
E_t^{(i_0,j_0,s_0)})_{t=1}^{L_j}$ by Lemma~2.7 to satisfy (2.14)--(2.17). 

For $i,j\le n$ define 
$$x(i,j) = {1\over k_0} \sum_{s=1}^{k_0} x(i,j,s)\ . 
\leqno(2.18)$$ 
Observe that for any $j\le n$, the family 
$(m_t^{(i,j,s)} ,E_t^{(i,j,s)})_{i\le n,\, s\le k_0,\, t\le L_j}$ is 
admissible (when ordered lexicographically, $(i,j,s,t)$) 
and so for $(a_i)_1^n\subseteq 
\real$, 
$$\eqalign{\Big\| \sum_{i=1}^n a_i \sum_{s=1}^{k_0} x(i,j,s)\Big\| 
& \ge \Big\| \sum_{i=1}^n a_i \sum_{s=1}^{k_0} x(i,j,s)\Big\|_{L_j\cdot n
\cdot k_0}\cr 
&\ge {1\over f(L_jnk_0)} \sum_{i=1}^n a_i \sum_{s=1}^{k_0} \sum_{t=1}^{L_j} 
\trivert E_t^{(i,j,s)} x(i,j,s)\trivert_{m_t^{(i,j,s)}} \cr 
&\ge {f(L_j) \over  f(L_j\cdot n\cdot k_0)} 
\sum_{i=1}^n\ \sum_{s=1}^{k_0} \|a_i x(i,j,s)\|_{(L_j,m_0^{(i,j,s)})}
\text{ (by (2.16))}\cr 
&\ge (1-\delta)^2 k_0 \sum_{i=1}^n |a_i| \ \text{ (by (2.16) and (2.12))}\ .
\cr}
\leqno(2.19)$$ 

We conclude from (2.18) and (2.19) that 
for all $1\le j\le n$, $(x(i,j))_{i=1}^n$ is 
$(1-\delta)^{-2}$-equivalent to the unit vector basis of $\ell_1^n$. 
By our choice of $\delta$ (2.8) we deduce the left hand inequality of (2.7). 

To prove the right hand estimates let $(a(i,j))_{i,j\le n}\subseteq \real$ 
with $\|\sum_{i,j} a(i,j) x(i,j)\|=1$ 
and let $\ell\in\nat$. The argument of (2.19) yields that for fixed $(i,j)$, 
$$ \Big\| \sum_{s=1}^{k_0} x(i,j,s)\Big\| 
\ge k_0 (1-\delta)^2 $$ 
and so each $x(i,j)$ is an $\ell_1^{k_0}$-average with constant 
$(1-\delta)^2$. Thus if $\ell\le L_0$, 
$$\Big\| \sum_{i,j=1}^n a(i,j) x(i,j)\Big\|_\ell 
< \Big\| \sum_{i,j=1}^n a(i,j) x(i,j)\Big\|$$ 
from Remark 2.6 and (2.11). 

If $\ell > L_0$ then there is at most one $j_0\le n$ so that $\ell \in 
[L'_{j_0-1} +1, L'_{j_0}]$. Then 
$$\eqalign{ 
\Big\| \sum_{\scriptstyle i,j=1\atop\scriptstyle j\ne j_0}^n 
a(i,j) x(i,j)\Big\|_\ell 
& \le \left[ {2n^2\over f(\ell)} + n^2\delta\right] \max \{ |a(i,j)| : 
i,j\le n\} \text{ (by (2.15) and (2.17))}\cr
&\le 3n^2\delta \max \{ |a(i,j)| : i,j\le n\} \text{ (by (2.10) and } 
\ell\ge L_0)\ .\cr}$$ 
Thus 
$$\eqalignno{\|\sum_{i,j=1}^n a(i,j) x(i,j)\Big\| 
&\le \max_{j\le n} \sum_{i=1}^n |a(i,j)| + 3n^2\delta\max 
\{ |a(i,j)| : i,j\le n\} \cr 
&\le (1+\varep) \max_{j\le n} \sum_{i=1}^n |a(i,j)|\ \text{ (by (2.9)).} 
&\qed\cr}$$ 
\enddemo 

We turn now to the task of proving Lemmas 2.5 and 2.7. 
This requires several steps. 

\proclaim{Step 1} 
Let $(y_i)$ be a block basis of $(e_i)$. 
There exists $1\le p<\infty$ so that $\ell_p$ is block finitely 
represented in $(y_i)$. 
\endproclaim 

\demo{Proof} 
By Krivine's theorem $\ell_p$ is block finitely represented in $(y_i)$ 
for some $1\le p\le\infty$. Suppose $p=\infty$ (otherwise we are done). 
Let $\varep >0$ and choose a block basis $(w_i)$ of $(y_i)$ as follows. 
Let $w_1\in \text{span}(y_i)$ be an $\ell_\infty^2$-average with constant 
${1+\varep}$ and let $m_1=2$. 
If $w_1,\ldots,w_k$ are chosen, let $m_{k+1}\in\nat$ 
with $f(m_{k+1})\ge \sum_{j=1}^k |\text{supp}(w)|$ and let $w_{k+1}$ 
be an $\ell_\infty^{m_{k+1}}$-average of a normalized block basis 
in span$(y_i)$ with constant $(1+\varep)$ and, of course $w_k<w_{k+1}$. 
For any $\ell\in\nat$ and $(k_i)\subseteq\nat$, $(\text{supp}(w_{k_i}), 
m_{k_i})_{i=1}^\ell$ is admissible and so 
$$\Big\| \sum_{i=1}^\ell w_{k_i} \Big\| 
\ge {1\over f(\ell)} \sum_{i=1}^\ell \trivert w_{k_i}\trivert_{m_{k_i}} 
\ge {\ell\over f(\ell)}\ .$$
Also $\|w_{k_i}\| \le 1+\varep$ for all $i$. We conclude from James' 
well known argument \cite{J} that $\ell_1$ is block finitely representable 
in $(y_i)$.\qed
\enddemo 

\proclaim{Step 2}
Let $1\le p<\infty$, $k,m,\ell\in \nat$ and let $x\in X$ be an 
$\ell_p^k$-average with constant $2$. Then 
\vskip1pt 
\iitem{a)} $\trivert x\trivert_m \le 4m^{-1/p} + k^{-1/p}$ 
\vskip1pt 
\iitem{b)} $\|x\|_\ell \le {1\over f(\ell)} [C_p +2\ell k^{-1/p}]$ where 
$C_p = \sum_{i=0}^\infty 4\cdot 2^{-i/p} = 4(1-2^{-1/p})^{-1}$. 
\endproclaim 

\demo{Proof} 
Let $x= k^{-1/p} \sum_{i=1}^k x_i$ where $(x_i)_{i=1}^k$ is a normalized 
block basis of $(e_i)$ which is 2-equivalent to the unit vector basis of 
$\ell_p^k$. Choose $E_1< E_2<\cdots < E_m$ so that $E_i\subseteq\text{supp}(x)$ 
for $i\le m$,\footnotemark\footnotetext{It may not be possible to 
insure that $E_i\subseteq \text{supp}(x)$ (e.g., if $|\text{supp}(x)|<m$) 
but this case presents only notational troubles.} 
$\trivert x\trivert_m = {1\over m} \sum_{i=1}^m 
\|E_ix\|$ and choose an admissible family $(m_i,E_i^\ell)_{i=1}^\ell$ 
so that $E_i^\ell \subseteq\text{supp }x$ for $i\le \ell$ and 
$$\|x\|_\ell = {1\over f(\ell)} \sum_{i=1}^\ell \trivert E_i^\ell x\trivert_
{m_i}\ .$$ 
For $i\le k$, set $I_i = \{j\le m: E_j \subseteq \text{supp}(x_i)\}$ and 
$I_0 = \{j\le m: j\notin\bigcup_{i=1}^k I_i\}$. 
For $j\in I_0$, let $k_j = \min \{i\le k:\text{supp}(x_i)\cap E_j 
\ne \emptyset\}$ and $K_j = \max \{i\le k:\text{supp}(x_i)\cap E_j 
\ne \emptyset\}$. Note that $|I_0| \le \min (k,m)$. 
Similarly define $I_i^\ell, I_0^\ell, k_j^\ell$ and $K_j^\ell$ from the 
collection $(E_j^\ell)_{j=1}^\ell$ and note that $|I_0^\ell|\le \min (k,\ell)$.
Note also that $k_j <K_j\le k_{j'}$ for $j<j'\in I_0$ 
and $k_j^\ell < K_j^\ell \le k_{j'}^\ell$ for $j<j'\in I_0^\ell$. 

In order to prove a) we note first the following estimate. 
$$\eqalign{\sum_{j\in I_0} \|E_jx\| 
& \le k^{-1/p} \sum_{j\in I_0} \Big\|\sum_{i=k_j}^{K_j} x_i\Big\|\cr 
&\le 2k^{-1/p} \sum_{j\in I_0} (K_j-k_j+1)^{1/p} 
\le 4k^{-1/p} \sum_{j\in I_0} (K_j-k_j)^{1/p}\cr 
&\le 4k^{-1/p} \biggl[ \sum_{j\in I_0} (K_j-k_j)\biggr]^{1/p} 
|I_0|^{1/q}\quad \Bigl(\text{ by H\"older's inequality with }
\dfrac1p +\dfrac1q =1 \Bigr) \cr
&\le 4|I_0|^{1/q}\quad \Bigl( \text{ using }\ \sum_{j\in I_0} (K_j-k_j) 
\le k\Bigr) \cr 
&\le 4m^{1-1/p}\ .\cr}$$ 
Next let $i\le k$ with $I_i\ne\emptyset$. Then 
$$\eqalign{ 
\sum_{j\in I_i} \|E_jx\| & = |I_i| k^{-1/p} |I_i|^{-1} 
\sum_{j\in I_i} \|E_j x_i\|\cr 
&\le |I_i| k^{-1/p} \trivert x_i\trivert_{|I_i|}\cr
&\le |I_i| k^{-1/p}\quad (\text{since } \trivert x_i\trivert_{|I_i|} 
\le \|x_i\| =1)\ .\cr}$$ 

Combining these two estimates we obtain 
$$\eqalign{\trivert x\trivert_m &\le {1\over m} \biggl[ 4m^{1-1/p} + 
\sum_{i=1}^k |I_i| k^{-1/p}\biggr] \cr 
&\le 4m^{-1/p} + k^{-1/p}\ ,\text{ which proves a).}\cr}$$

Next we verify b).  
$$\eqalign{ \sum_{j\in I_0^\ell} \trivert E_j^\ell x\trivert_{m_j} 
&\le k^{-1/p} \sum_{j\in I_0^\ell} \trivert \sum_{i=k_j^\ell}^{K_j^\ell} 
x_i\trivert_{m_j}\cr 
&\le k^{-1/p} \sum_{j\in I_0^\ell} (K_j^\ell - k_j^\ell +1)^{1/p} 
\bigl[ 4m_j^{-1/p} + (K_j^\ell - k_j^\ell +1)^{-1/p}\bigr]\ \text{ (by a))}\cr 
&\le \sum_{j\in I_0^\ell} 4m_j^{-1/p} + |I_0^\ell | k^{-1/p}\cr 
&\le C_p + \ell k^{-1/p}\ .\cr}$$ 
Furthermore for $i\le k$, 
$$\sum_{j\in I_i^\ell} \trivert E_j^\ell x\trivert_{m_j} 
= k^{-1/p}\sum_{j\in I_i^\ell} \trivert E_jx_i\trivert_{m_j} 
\le k^{-1/p} |I_i^\ell|\ .$$ 
Since $\sum_{i=1}^k |I_i^\ell| \le \ell$ we deduce that 
$$\sum_{j=1}^\ell \trivert E_j^\ell x\trivert_{m_j} \le C_p +2\ell k^{-1/p}$$ 
which yields b).\qed 
\enddemo 

Lemma 2.5 follows directly from the next step. 

\proclaim{Step 3} 
Let $n\in \nat$, $1\le p<\infty$ and let $(x_i)_{i=1}^n$ be a block basis 
of $(e_i)$ consisting of $\ell_p^{k_i}$-averages with constant $2$ for some 
$(k_i)_1^n\subseteq \nat$. Let $(m_i,E_i)_{i=1}^\ell$ be an admissible family 
and $(a_i)_{i=1}^n \subseteq [-1,1]$. 
Let $k_0=\min \{k_i:i\le n\}$ and let $j_0\le \ell$ be maximal so that 
$\sum_{j=1}^{j_0-1} |E_j| \le k_0^{1/2p}$. 
Then for $x= \sum_{i=1}^n a_i x_i$, 
$$\sum_{\scriptstyle j=1\atop\scriptstyle j\ne j_0}^\ell \trivert E_j x
\trivert_{m_j} \le 6n\ell k_0^{-1/2p}\ .$$ 
\endproclaim 

\demo{Proof}
Since $\|x\|_\infty \le k_0^{-1/p}$, by our choice of $j_0$, 
$$\sum_{j=1}^{j_0-1}\trivert E_jx\trivert_{m_j}
\le k_0^{1/2p} k_0^{-1/p} = k_0^{-1/2p}\ .$$
If $j_0<\ell$ then 
$$k_0^{1/2p} \le \sum_{j=1}^{j_0} |E_j| < m_{j_0+1} \le m_j
\ \text{ for all }\ j\ge j_0+1\ .$$ 
>From Step 2 a) and this we obtain 
$$\eqalign{\sum_{j=j_0+1}^\ell \trivert E_jx\trivert_{m_j} 
&\le \sum_{i=1}^n \sum_{j=j_0+1}^\ell \trivert x_i\trivert_{m_j}\cr 
&\le n\ell [k_0^{-1/p} + 4k_0^{-1/2p}] \cr 
&\le 5n\ell k_0^{-1/2p} \ .\cr}$$
Step 3 follows immediately from these two estimates.\qed
\enddemo 

The next step along with Step 1 yields that $c_0$ is block finitely 
represented in every block basis of $(e_i)$ and hence, from the proof of 
Step~1, so is $\ell_1$. 

\proclaim{Step 4} 
Let $1\le p<\infty$, $0<\varep<{f(2)-1\over 2}$, 
$n\in\nat$ and let $(y_i)_{i=1}^n $ be 
a block basis of $(e_i)$ consisting of $\ell_p^{k_i}$-averages with constant 
$1+\varep$. Suppose, in addition, that 
$$f\left( {\varep k_1^{1/2p} \over 6n}\right) \ge 
{n(C_p+2)\over \varep} 
\leqno(2.20)$$ 
and for $2\le i\le n$, 
$$f(k_i) > {p\over\varep} \sum_{s=1}^{i-1} |\text{\rm supp}(y_s)|\ .
\leqno(2.21)$$ 
Let $y=\sum_{i=1}^n y_i$. Then 
$$\|y\|_\ell \le {1\over f(\ell)} \bigl[\|y\|+\varep\bigr] 
\ \text{ if }\ \ell \le {\varep k_1^{1/2p} \over 6n} 
\leqno\text{\rm a)}$$
and 
$$\|y\|_\ell \le 2\varep +\max_{i\le n} \|y_i\|_\ell\ \text{ if }\ 
\ell > {\varep k_1^{1/2p}\over 6n} 
\leqno\text{\rm b)}$$ 
In particular $\|y\| \le 2\varep +\max_{i\le n} \|y_i\|$.
\endproclaim 

\demo{Proof} 
a) follows immediately from Step 3. Now suppose $\ell > {\varep k_1^{1/2p} 
\over 6n}$. 
Let $(m_j,E_j)_{j=1}^\ell$ be an admissible family so that 
$$\|y\|_\ell = {1\over f(\ell)} \sum_{j=1}^\ell \trivert E_j\trivert_{m_j}\ .$$ 
Let $i_0\ge 1$ be maximal so that 
$$\sum_{s=1}^{i_0-1} |\text{supp}(y_j)| \le \varep\, f(\ell)\ .$$ 
Note that 
$${1\over f(\ell)} \sum_{j=1}^\ell \trivert E_j\biggl( \sum_{i=1}^{i_0-1} 
y_i\biggr)\trivert_{m_j} \le {1\over f(\ell)} \sum_{i=1}^{i_0-1} 
|\text{supp}(y_i)| \le \varep\ . 
\leqno(2.22)$$ 
Note also that if $i_0<n$ then by (2.21),  
$$\varep \,f(\ell) \le \sum_{s=1}^{i_0} |\text{supp}(y_s)| 
\le {\varep\over p} f(k_{i_0+1})$$ 
hence $\ell \le k_{i_0+1}^{1/p}$ and so $\ell k_i^{-1/p} < 1$ 
for $i>i_0$. 
Now from Step~2~b) we have 
$$\eqalign{
{1\over f(\ell)} \sum_{j=1}^\ell \trivert E_j\biggl( \sum_{i=i_0+1}^n 
y_i\biggr) \trivert_{m_j} &\le \sum_{i=i_0+1}^n \|y_i\|_\ell \cr 
&\le {1\over f(\ell)} \sum_{i=i_0+1}^n (C_p +2\ell \, k_i^{-1/p}) 
\le {n(C_p+2)\over f(\ell)} \cr 
&\le {n(C_p+2) \over f(\varep\, k_1^{1/2p} (6n)^{-1})} 
\le \varep\ ,\ \text{ (by (2.20)).}\cr} 
\leqno(2.23)$$ 
Hence from (2.22) and (2.23), 
$$\eqalign{\|y\|_\ell & \le\varep + {1\over f(\ell)} 
\sum_{j=1}^\ell \trivert E_j y_{i_0}\trivert_{m_j} +\varep\cr 
&\le 2\varep + \|y_{i_0}\|_\ell\ .\cr}$$ 

To see the ``in particular'' statement we note that if $\|y\|=\|y\|_\ell$ 
for some $\ell\in \nat$ then $\|y\|=\|y\|_\ell$ for some $\ell\ge2$ by 
Remark~2.4. If $2\le \bar \ell\le {\varep k_1^{1/2p}\over 6n}$ then 
since $\|y\| >\tfrac12$, a) and $\varep < {f(2)-1\over 2}$ yields  that 
$\|y\|_{\bar\ell} < \|y\|$. 
Thus b) yields the assertion.\qed
\enddemo 

We have one final step before proving Lemma 2.7. 

\proclaim{Step 5} 
Let $m\in\nat$ and $\varep>0$. There exists $\delta = \delta (m,\varep)>0$ 
so that whenever $(z_i)_{i=1}^m$ is a block basis of $(e_i)$ satisfying 
\vskip1pt
{\rm (A)} For each $i\le m$, $z_i= \sum_{j=1}^{n_i} z(i,j)$, 
where $(z(i,j))_{j=1}^{n_i}$ is a block basis of $(e_i)$ consisting of 
$\ell_1^{k(i,j)}$-averages with constant $1+\delta$ satisfying for $i\le m$, 
$$f\left( {\delta k(i,1)^{1/2p}) \over 6n_i}\right) 
> {n_i(C_1+2)\over\delta} 
\leqno(2.20)'$$ 
and 
$$f\bigl( k(i,j)\bigr) > {1\over\delta} \sum_{s=1}^{j-1} |\text{\rm supp} 
\bigl(z(i,s)\bigr)|\ ,\quad \text{for }\ j\ge 2 
\leqno(2.21)'$$ 
\vskip1pt
{\rm (B)} $n_1>{m\over\delta}$ and $f(n_i)>\sum_{j=1}^{i-1} |\text{\rm supp } 
z_j|$ for $2\le i\le m$. 
\vskip1pt 
Then for all $\ell \in\nat$, 
$$\Big\| \sum_{i=1}^m z_i\Big\|_\ell \le (1+\varep) \max 
\left\{ 1,{m\over f(\ell) f(m/\min (\ell,m)) }\right\}$$ 
and thus (since $f(xy) \le f(x)f(y)$ for $x,y\ge1$, 
see e.g., \cite{S1} lemma 1) 
$$\Big\|\sum_{i=1}^m z_i\Big\|_\ell \le (1+\varep) {m\over f(m)}\ .$$
\endproclaim 

\demo{Proof of Lemma 2.7} 

Let $\tfrac13 >\varep >0$, $L_0\in\nat$. Choose $L_1>L_0$ with 
${f(L_1)\over f(L_1/L_0)} \le 1+\varep$ and ${f(L_1)\over L_1}<\varep/2$. 
Then take $L'_1 >L_1$ with 
${1\over 1-\varep}({f(L_1)\over  L_1} + {f(L_1)\over f(L'_1)}) <\varep$. 

Let $m_0\in\nat$ and let $(y_i)$ be a block basis of $(e_i)$. 
Let $(z_i)_{i=1}^{L_1}$ be a block basis of $(y_i)$ which satisfies 
A) and B) of Step~5 for ${\varep\over L_1} \equiv \varep'$ and 
$\delta = \delta (L_1,\varep')$. We may also assume that $n_1\ge m_0$. 
This can be done since $\ell_1$ is block finitely represented in $(y_i)$. 

>From Step 5 and the properties of $f$,  
$$\Big\|\sum_{i=1}^{L_1} z_i\Big\| 
\le (1+\varep') \sup_{\ell\in\nat} {L_1\over  f(\ell) f(L_1/\min(L_1,\ell) )} 
\le {(1+\varep') L_1\over f(L_1)}
\leqno(2.24)$$ 

Set 
$$\bar x= {1\over 1+\varep'} \ {f(L_1)\over L_1} 
\sum_{i=1}^{L_1} z_i \in Ba(X)\ \text{ and }\ 
x= \bar x\slash \|\bar x\|\ .$$ 
Using the notation of Step 5, 
$$\eqalign{1\ge \|\bar x\| & \ge \|\bar x\|_{(L_1,m_0)}\cr 
&\ge {1\over 1+\varep'} \ {f(L_1)\over L_1} \ {1\over f(L_1)} 
\sum_{i=1}^{L_1} \trivert z_i\trivert_{n_i}\cr 
&\ge {1\over 1+\varep'}\ {1\over L_1} \sum_{i=1}^{L_1} {1\over n_i} 
\sum_{j=1}^{n_i} \|z(i,j)\|\cr 
&\ge {1\over 1+\varep'}\ {1\over L_1} \sum_{i=1}^{L_1} {1\over n_i} 
\sum_{j=1}^{n_i} {1\over 1+\varep'} 
= {1\over (1+\varep')^2} > 1-\varep\ .\cr}$$ 
Thus $\|x\|_{(L_1,m_0)} > 1-\varep$ which proves b) of Lemma 2.7. 

For $\ell\le L_0$, 
$$\eqalign{ \|\bar x\|_\ell & \le {1\over 1+\varep'} \ {f(L_1)\over L_1}\ 
{L_1(1+\varep') \over f(\ell)f({L_1\over \ell}) }\ \text{ (by (2.24))}\cr
&= {f(L_1)\over f(\ell) f({L_1\over\ell})} 
\le {f(L_1)\over f(\ell) f({L_1\over L_0})} 
\le {1+\varep \over f(\ell)}\ \text{ and so}\cr 
\|x\|_\ell & < {1+\varep\over f(\ell)}\ {1\over 1-\varep} < {2\over f(\ell)} 
\ \text{ which proves a).}\cr}$$ 

Finally if $\ell \ge L'_1$, 
$$\eqalign{\|\bar x\|_\ell & \le {1\over 1+\varep'} \ {f(L_1)\over L_1} 
(1+\varep') \max \left\{ 1,{L_1\over f(\ell)}\right\}\ \text{ (from Step 5)}\cr 
&\le {f(L_1)\over L_1} + {f(L_1)\over f(\ell)} 
\le {f(L_1)\over L_1} + {f(L_1)\over f(L'_1)} <\varep \cr}$$ 
and so 
$$\|x\|_\ell <{1\over 1-\varep}\left( {f(L_1)\over L_1} + 
{f(L_1)\over f(L'_1)}\right) <\varep$$ 
which proves c).\qed 
\enddemo 

\demo{Proof of Step 5} 

We proceed by induction on $m$. Let ${f(2)-1\over 2}>\varep >0$. 
For $m=1$ the conclusion we desire becomes 
$$\|z_1\|_\ell \le 1+\varep \text{ which holds by Step 4 if } 
\delta < {\varep\over 2}\ .$$ 

Assume Step 5 has been proved for all $1\le m'<m$. Let $\delta>0$ be fixed 
small enough to satisfy the conclusion of Step~5 for all $m'<m$ and an 
$\varep'<\varep$ to be specified later. 
and let $(z_i)_{i=1}^m$ satisfy the hypothesis of Step~5. 
\smallskip

\noindent {\it Case 1\/}: 
$\ell \ge \dfrac{\delta k(1,1)^{1/2}}{6n_1}$. 

Our growth conditions imply that 
$$\left({\delta k(i,1)^{1/2} \over 6n_i}\right)$$ 
is an increasing sequence. Choose $i_0\le m$ to be maximal so that 
$${\delta k(i_0,1)^{1/2}\over 6n_{i_0}} \le \ell\ .$$ 
We have that 
$$\eqalign{ f(\ell) & \ge f\left( {\delta k(i_0,1)^{1/2} \over 6n_{i_0}} 
\right) \ge {n_{i_0}\over \delta} \ \text{ (from  (2.20)$'$)}\cr 
&\ge \delta^{-1} \sum_{i=1}^{i_0-1} |\text{supp}(z_i)|\ \text{ (by B)}\cr}$$ 
which implies that 
$$\Big\| \sum_{i=1}^{i_0-1} z_i\Big\|_\ell 
\le {1\over f(\ell)} \sum_{i=1}^{i_0-1} |\text{supp}(z_i)| \le \delta\ . 
\leqno(2.25)$$ 

Furthermore, 
$$\eqalign{\Big\|\sum_{i=i_0+1}^m z_i\Big\|_\ell 
&\le {1\over f(\ell)} \sum_{i=i_0+1}^m \bigl( \|z_i\|+\delta\bigr)\cr 
&\qquad \text{(from Step 4 a) since }  \ell \le {\delta k(i,1)^{1/2} 
\over 6n_i}\ \text{ for }\ i>i_0)\cr
&\le {m\over f(\ell)} (1+3\delta) \le {m(1+3\delta) \over 
f({\delta k(1,1)^{1/2}\over 6n_1} )}\ \text{ (using $\|z_i\|< 1+2\delta$ 
from Step 4b)}\cr 
&\le {m(1+3\delta)\over n_1(C_1+2)} \delta\ \text{ (from (2.20)$'$)} 
\le \left( {1+3\delta\over C_1 +2}\right) \delta \cr 
&\le \delta\ \text{ (provided }\ \delta \le 1/2)\ .\cr}$$ 

>From this and (2.25) we have 
$$\Big\|\sum_{i=1}^m z_i\Big\|_\ell \le 2\delta + \|z_{i_0}\|_\ell 
\le 1+  4\delta\ . 
\leqno(2.26)$$ 

\noindent {\it Case 2\/}: $\ell \le \dfrac{\delta k(1,1)^{1/2}}{6n_1}$. 

Let $(m_j,E_j)_{j=1}^\ell$ be an admissible family 
for which if $z=\sum_{i=1}^m z_i$, then 
$$\|z\|_\ell = {1\over f(\ell)} \sum_{j=1}^\ell \trivert E_j z\trivert_{m_j}
\ .$$ 
For $1\le i\le m$ let $j_i \le \ell$ be maximal so that 
$\sum_{j=1}^{j_i-1} |E_j| \le k(i,1)^{1/2}$. 

It is possible that $j_s=j_t$ for some $t\ne s$. Let $\ell' \equiv 
|\{j_1,\ldots,j_m\}|$ and relabel the set $\{E_{j_1},E_{j_2},\ldots, 
E_{j_m}\}$ as $E'_1 <E'_2 < \cdots < E'_{\ell'}$. 
Choose intervals $I_i\subseteq \nat$ so that $\bigcup_1^m I_i = [0,\max 
\text{ supp}(z)]$ and $\text{supp}(z_i) \subseteq I_i$ for $i\le m$. 
Define for $s\le \ell'$, 
$$\widetilde E_s \equiv E'_s\setminus \cup \{ I_i : E_{j_i} \ne E'_s\}\ .$$ 
Note that for all $i\le m$, $\widetilde E_s\cap I_i\ne\emptyset$ 
for at most one $s$. Thus we can choose $0=k_0<k_1<\cdots <k_{\ell'}$ 
so that $\widetilde E_s\subseteq \bigcup_{i=k_{s-1}+1}^{k_s} I_i$ 
for $s\le \ell'$. 
Note also that for $j\in\{j_1,\ldots,j_m\}$, if $E'_s =E_j$ then 
$$\eqalign{ \trivert E_jz\trivert_{m_j} 
&\le \trivert\widetilde E_s z\trivert_{m_j} 
+ \sum_{\scriptstyle i\le m\atop\scriptstyle E_{j_i}\ne E'_s}  
\trivert E_jz_i\trivert_{m_j}\cr 
&\le \|\widetilde E_s z\| + \sum_{\scriptstyle i\le m\atop\scriptstyle 
E_{j_i}\ne E'_s} \trivert E_j z_i\trivert_{m_j}\ .\cr}$$ 
This implies that 
$$\eqalign{ \sum_{j=1}^\ell \trivert E_jz\trivert_{m_j} 
& = \sum_{j\in \{j_1,\ldots,j_m\}} \trivert E_jz\trivert_{m_j} 
+ \sum_{\scriptstyle j\le \ell\atop\scriptstyle j\notin \{j_1,\ldots,j_m\}} 
\trivert E_j z\trivert_{m_j}\cr 
&\le \sum_{s=1}^{\ell'} \|\widetilde E_s z\| +\sum_{j\in \{j_1,\ldots,j_m\}} 
\sum_{\scriptstyle i\le m\atop\scriptstyle E_{j_i}\ne E_j} 
\trivert E_j z_i\trivert_{m_j} 
+ \sum_{\scriptstyle j\le \ell\atop\scriptstyle j\notin \{j_1,\ldots,j_m\}} 
\trivert E_j z\trivert_{m_j}\cr 
&\le \sum_{s=1}^{\ell'} \|\widetilde E_s z\| 
+ \sum_{i=1}^m \biggl( \sum_{\scriptstyle j\in \{j_1,\ldots,j_m\}\atop 
\scriptstyle E_{j_i}\ne E_j} \trivert E_j z_i\trivert_{m_j} 
+ \sum_{\scriptstyle j\le \ell\atop\scriptstyle j\notin \{j_1,\ldots,j_m\}} 
\trivert E_j z_i\trivert_{m_j}\biggr) \cr 
&\le \sum_{s=1}^{\ell'} \|\widetilde E_s z\| + \sum_{i=1}^m 
\sum_{\scriptstyle j\le \ell\atop\scriptstyle j\ne j_i} 
\trivert E_j z_i\trivert_{m_j}\cr 
&\le \sum_{s=1}^{\ell'} \Big\| \sum_{i=k_{s-1}+1}^{k_s} z_i \Big\| 
+ \sum_{i=1}^m 6n_i\ell\, k(i,1)^{-1/2}\text{ (by Step 3)}\cr
&\le \sum_{j=1}^{\ell'} \Big\|\sum_{i=k_{j-1}+1}^{k_j} z_i\Big\| 
+m\delta\cr
&\qquad \text{(using that $6n_i\ell\,k(i,1)^{-1/2}\le 6n_i\ell\,k(1,1)^{-1/2} 
<\delta$ by Case 2).}\cr}$$ 
If $\ell' = 1$ we deduce that 
$$\|z\|_\ell \le {\|z\| +m\delta \over f(\ell)}\ .
\leqno(2.27)$$ 
If $\ell' >1$ since $\delta$ has been chosen smaller than $\delta (m',\varep')$ 
for all $m' <m$, we deduce that 
$$\eqalign{
\sum_{j=1}^{\ell'} \Big\| \sum_{i=k_{j-1}+1}^{k_j} z_i\Big\| 
&\le (1+\varep') \sum_{j=1}^{\ell'} {k_j-k_{j-1} \over f(k_j-k_{j-1})} \cr
&\le (1+\varep') {m\over f({m\over \ell'})} 
\le (1+\varep') {m\over f(m/\min (m,\ell))}\ .\cr}
\leqno(2.28)$$ 
The next to the last inequality above derives from the fact that the maximum 
of $g(a_1,\ldots, a_{\ell'})$ $\equiv \sum_{i=1}^{\ell'} {a_i\over f(a_i)}$ 
on the set $\{(a_i)_1^{\ell'} : a_i \ge 1$ and $\sum_1^{\ell'} a_i=m\}$ 
is achieved for $a_1=\cdots=a_{\ell'} = {m\over \ell'}$. 

>From (2.26), (2.27) and (2.28) we have 
$$\left\{\eqalign{\|z\| &= \sup_\ell \|z\|_\ell 
\le \max \left\{ 1+4\delta, \sup_{\ell\in\nat} {1+\varep'\over f(\ell)} 
\left[ {m\over f(m/\min (m,\ell))} +m\delta\right]\right\}\cr 
&\le (1+\varep') \left[ {m\over f(m)} +m\delta\right]\ .\cr} \right. 
\leqno(2.29)$$ 
Thus, again using (2.26)--(2.28) and (2.29), and assuming $\varep' <\tfrac12$, 
$$\eqalign{\|z\|_\ell 
&\le \max \left\{ 1+4\delta, {\|z\|+m\delta\over f(\ell)} , 
{1+\varep'\over f(\ell)} \left[ {m\over f(m/\min (m,\ell))} 
+m\delta\right]\right\} \cr 
&\le 4m\delta +\max \left\{ 1+4\delta ,(1+\varep') {m\over f(m)f(\ell)}, 
(1+\varep') {m\over f(m/\min (m,\ell))} \right\}\ .\cr}$$ 
Step 5 follows by taking $\varep'$ and $\delta$ sufficient 
small.\qed 
\enddemo 

\demo{Remark 2.8}
The authors have recently obtained a conditional version of the example in 
Section~2: a basic sequence $(x_i)$ with the property that if $(y_i)$ is 
any block basis of $(x_i)$, $\varep>0$ and $(z_i)_1^n$ is a monotone 
basic sequence then there exists a block basis $(w_i)_1^n$ of $(y_i)$ 
which is $(1+\varep)$-equivalent to $(z_i)_1^n$. 
\enddemo

\enddocument 

\end